\documentclass[12pt,leqno]{article}
\usepackage{indentfirst}

\usepackage{amsmath}
\usepackage{amscd}
\usepackage{amssymb}
\usepackage{amsthm}

\def\zR{\ensuremath{\mathbb{R}}}

\newcommand{\Proof}{{\it Proof}}

\numberwithin{equation}{section}
\renewcommand{\theequation}{\thesection.\arabic{equation}}
\makeatletter
\newtheoremstyle{estiloobs}
  {6pt}
  {0pt}
  {\upshape}
  {}
  {}
  {}
  {.5em}
  {(\thmnumber{#2}) \thmname{\textbf{#1}}
   \thmnote{({#3})}}
\makeatother
\theoremstyle{estiloobs}
\newtheorem{remark}[equation]{Remark:}

\makeatletter
\newtheoremstyle{estilotheorem}
  {6pt}
  {0pt}
  {\slshape}
  {}
  {}
  {}
  {.5em}
  {(\thmnumber{#2}) \thmname{\textbf{#1}}
   \thmnote{({#3})}}
\makeatother
\theoremstyle{estilotheorem}

\newtheorem{lemma}[equation]{Lemma:}
\newtheorem{corollary}[equation]{Corollary:}
\newtheorem{theorem}[equation]{Theorem:}
\makeatletter
\newtheoremstyle{estilolista}
  {0pt}
  {0pt}
  {\slshape}
  {}
  {}
  {}
  {.5em}
  {(\thmnumber{#2})\thmname{\textbf{#1}}
   \thmnote{({#3})}}
\makeatother
\theoremstyle{estilolista}

\newcommand{\pedro}{\theequation}

\begin{document}

\title{Boundedness of fractional operators in weighted variable exponent spaces
with non doubling measures}

\vskip 0.3 truecm

\author{Osvaldo Gorosito-Gladis Pradolini-Oscar Salinas \thanks{The authors were
supported by Consejo Nacional de Investigaciones Cient\'\i ficas y
T\'ecnicas de la Rep\'ublica Argentina and Universidad Nacional
del Litoral.\newline \indent Keywords and phrases: variable
exponent, weighted spaces, non doubling measures
\newline \indent 1991 Mathematics Subject Classification: Primary
42B25.\newline }}

\date{\vspace{-1.5cm}}

\maketitle

\begin{abstract}
In the context of variable exponent Lebesgue spaces equipped with a lower Ahlfors measure we obtain
  weighted norm inequalities over bounded domains for the
centered fractional maximal function
 and the fractional integral operator.

\end{abstract}
\section{Introduction and statement of the main results}

During the last two decades the variable exponent Lebesgue spaces
$L^{p(x)}$ have been studied intensively. They seem to be the most
adequate context for studying a great variety of problems related
to certain classes of fluids that are characterized by their
ability to undergo significant changes in their mechanical
properties when an electric field is applied (see \cite{R}).

Recently, a number of authors, interested in studying the
continuity of certain classical operators from harmonic analysis,
have succeeded in proving some boundedness results in such spaces.
In fact, un-weighted strong inequalities for the Hardy-Littlewood
and the fractional maximal operators were obtained in the
euclidean setting, under certain property of regularity on the
exponent. Particularly, Lars Diening proved the continuity of the
Hardy-Littlewood maximal operator in $\zR^n$ by requiring an
additional property of constancy on the exponent outside a fixed
ball (see \cite{D}). In proving his result, the technique
developed by the author differs from the one applied in the
classical theory, essentially based on interpolation. Cruz Uribe,
Fiorenza and Neugebauer took advantage of these techniques and
improved Diening's result by additionally assuming some type of
logarithm decay on the exponent. Thus they obtained norm
inequalities over open subsets of $\zR^n$ (see \cite{CUFN}).

Under the hypothesis of continuity and logarithm decay on the
exponent, the well known strong type inequality for the fractional
maximal operator was proved in the variable context over open
subsets of $\zR^n$ (see \cite{CCUF}).

In addition, boundedness results for the Hardy-Littlewood maximal
operator were obtained by Harjulehto, H\"{a}st\"{o} and Latvala in
$\zR^n$ over bounded domains with the novelty of using a non
necessarily doubling measure (see \cite{HHL}).

On the other hand, for the maximal operator, weighted norm
inequalities involving power weights, were obtained by
Kokilashvili and Samko
(\cite{KS}) in $\zR^n$ over open bounded domains.

In this article we prove weighted strong inequalities for
fractional operators. For maximal operators we include the case of the classical
Hardy-Littlewood maximal function in the setting of variable
exponent spaces defined over  bounded subsets of $\zR^n$ which have been equipped with
 a non necessarily doubling measure.  Such inequalities
provide a weighted version of those contained in \cite{HHL} for
the Hardy-Littlewood maximal function as well as those proved in
\cite{CCUF} for the fractional maximal in the standard context of
Lebesgue measurable spaces.

Additionally
 the class of weights involved in our results is wider
than that of power functions considered in \cite{KS} for the case
of the Hardy-Littlewood maximal operator referred to the
Lebesgue's measure. It is worth mentioning that a weighted pointwise relationship between the fractional and
 the Hardy-Littlewood maximal functions is also proved and it is not only interesting in itself but
 essential for the proof of one of  our main results as well.

 In the same context as the one described at the beginning, i.e.,
over generalized Lebesgue spaces equipped with a non-doubling measure, we also obtain
 weighted strong inequalities for the integral fractional operator. A version of
 Welland's
inequality in this variable setting is also given and it turns out
to play a fundamental role in proving the boundedness of that
operator.

 In the context of Lebesgue standard measure spaces we
give, in some sense, certain type of reverse H{\"{o}}lder
inequality which proves to be appropriate to obtain a special
class of weights for which the continuity of the fractional
integral holds. We want to point out that this class is larger
than that of power weights, generalizing in this way the results
due to Samko in \cite{S}. Finally we should say that the
techniques developed by the author, essentially based on a Hedberg
type inequality, differs
 from ours since Welland's inequality allows us to prove our results.

\bigskip

We first introduce the context in which we shall develop our
results.\\

Throughout this paper $Q=Q(x,l(Q))$ will denote a cube centered at
$x$ with length side $l(Q)$ and whose sides are parallel to the
coordinate axes. Moreover, $C$ will denote a positive constant not
necessarily the same on each occurrence.

Let us now consider a non- negative Borel regular measure $\mu$
defined over subsets in $\zR^n$. If $\Omega$ is a bounded
$\mu$-measurable set and $\beta:\Omega\rightarrow (0,\infty)$ is
a bounded function, we shall say that $\mu$ is a lower Ahlfors
$\beta(\cdot)$-regular measure in $\Omega$ if there exists a
positive constant $c$ such that the inequality
\begin{equation}\label{la}
c\leq\frac{\mu(Q(x, l(Q)))}{l(Q)^{\beta(x)}}\ ,
\end{equation}
holds for every $x\in \Omega$ and for every cube $Q\subset \zR^n$
such that $0<l(Q)<diam(\Omega)$. Particularly, whenever $\beta$
is a constant function we shall simply say that $\mu$ is lower
Ahlfors $\beta$-regular in $\Omega$.

It is clear that lower Ahlfors $\beta(\cdot)$-regular measures
have not necessarily got the doubling property, which means that
the inequality $\mu(2Q)\le C\mu(Q)$ holds for some constant $C\ge
1$ and for every cube $Q$.

The interest in studying these measures appears in connection with
the notion of the dimension of a metric space. By dimension it is
understood some quantity relating the measure of a cube with its
length side. Examples of measures with variable dimensions are
given in \cite{HHL}.

It is easy to see that lower Ahlfors $\beta(\cdot)$-regularity
implies lower Ahlfors $\beta^*$-regularity, where
$\beta^*=\sup_{\Omega}\beta$. In this paper we shall consider
lower Ahlfors $\beta$-regularity. It is not difficult to prove
that all our results imply the
corresponding results for the lower $\beta(\cdot)$-regularity case. In this article, $\beta$ will denote a
positive number related to the operators we shall be working with.

\bigskip

We now introduce the functional space we are going to deal
with. For additional information see \cite{KR}.\\

A non-negative $\mu$-measurable function $p$ from $\Omega$ to
$[1,\infty)$ is called an exponent. For simplicity we write
$p_*=\inf_\Omega p$ and $p^*=\sup_\Omega p$. An exponent is
bounded if $p^*<\infty$. We shall also assume that $p_*>1$. For
any set $A\subset \Omega$ we denote $(p_A)_*=\inf_{A} p$ and
$(p_A)^*=\sup_{A} p$.

For any $\mu$-measurable function $f$, the modular
$\rho_{p(\cdot),\Omega}$ is defined by
\begin{equation*}
\rho_{p(\cdot),\Omega}(f)=\int_{\Omega}|f(x)|^{p(x)}\, d\mu(x).
\end{equation*}
and the formula
\begin{equation*}
\|f\|_{p(.),\Omega}=\inf\{\lambda>0:
\rho_{p(\cdot),\Omega}(f/\lambda)\leq 1\}
\end{equation*}
is seen to define a norm.

The variable exponent Lebesgue space $L^{p(\cdot)}(\Omega)$
consists of those $\mu$-measurable functions $f$ supported in
$\Omega$ for which $\|f\|_{p(.),\Omega}<\infty$. If $p$ is an
exponent such that $p_*>1$, let $p'$ be the function defined by
$1/p(x)+1/p'(x)=1$. Topics related to general properties of this
space are treated in \cite{KR}. In particular, the generalized
H{\"{o}}lder inequality (see \cite{KR})
\begin{equation}\label{holder}
\int_{\Omega}|fg|\, d\mu \leq C
\|f\|_{p(\cdot),\Omega}\|g\|_{p'(\cdot),\Omega}
\end{equation}
holds and it shall be useful in our proofs.

 We shall deal with a class of bounded exponents which satisfy
certain property of regularity stronger than uniform continuity.
More precisely, an exponent $p$ is said to be log-H\"{o}lder
continuous if it satisfies the following inequality
\begin{equation*}
|p(x)-p(y)|\leq \frac{C}{\log(1/|x-y|)}\ ,\quad  x,\ y\in \Omega,
\quad |x-y|\leq 1/2.
\end{equation*}

It is worth mentioning that this condition guarantees regularity
results on variable exponent spaces. In \cite{D}, the author
proves that this condition along with the additional assumption
that $p$ is constant outside a fixed ball is sufficient for the
maximal operator to be bounded in $L^{p(\cdot)}(\zR^n)$. Moreover
in \cite{PR} it is shown that the boundedness of this operator
might fail for a general exponent $p$. In fact, the authors
proved that the modulus of continuity is optimal.

\vspace{0.5cm}

We finally introduce the maximal functions we are interested in
along with the class of weights involved with their properties of
boundedness.

Let $\mu$ be a lower Ahlfors $\beta$-regular measure. For $0\leq
\alpha<\beta$ the centered fractional maximal of a locally
integrable function $f$ is defined by
\begin{equation}\label{maxfrac}
M_{\alpha}f(x)=\sup_{r>0}\frac{1}{\mu(Q(x,r))^{1-\alpha/\beta}}\int_{Q(x,r)}
|f(y)|\, d\mu(y)
\end{equation}
where $Q(x,r)$ denotes a cube centered at $x$ with side-length
equal to r. If $\alpha=0$ in ($\ref{maxfrac}$) we simply write
$M_0=M$ for the classical Hardy-Littlewood maximal function.\\

A version of the fractional integral operator associated to the
maximal defined above is given by
\begin{equation}\label{intfrac}
I_{\alpha}f(x)=\int_{\Omega}
\frac{f(y)|x-y|^{\alpha}}{\mu(Q(x,2|x-y|))}\, d\mu(y).
\end{equation}

Let  $s$ be a real number such that $1<s< \infty$. We say that a
weight $w$ belongs to the $A_{s}(\Omega)$ class if there exists a
positive constant $C$ such that the inequality

\begin{equation*}
\left( \frac{1}{\mu(Q)}\int_{Q\cap \Omega}w\, d\mu \right) \left(
\frac{1}{\mu(Q)}\int_{Q\cap \Omega}w^{-1/(s-1)}\, d\mu
\right)^{s-1}\leq C
\end{equation*}
holds for each cube $Q$ centered at a point in $\Omega$.

If $\mu$ is the classical Lebesgue measure in $\zR^n$,
$0<\alpha<n$, $1<p<n/\alpha$, $1/q=1/p-\alpha/n$ and  $s=1+q/p'$,
it is well known that the $A_{s}(\zR^n)$ class characterizes the
boundedness of $M_{\alpha}$ from $L^p(w^{p/q})$ into $L^q(w)$ (see
for example \cite{MW}). Particularly if $\alpha=0$ and
$1<p<\infty$ this result gives the
boundedness of $M$ in $L^p(w)$.\\

Before stating our main results we introduce some additional
notation.

Given a continuous function $t$ defined in $\Omega$ we shall
denote
\begin{equation*}
\Omega_r^t=\{x\in \Omega : t(x)>r\}
\end{equation*}
for each $r\in\zR$. Let us observe that this set is not empty
whenever $ r< t^*$. Given $\epsilon>0$, related to this set, we
define
\begin{equation}\label{omega}\Omega_{r,\epsilon}^t=\Omega_r^t - \bigcup_{ x\in
\Omega-\Omega_r^t}B(x,\epsilon).
\end{equation}
It is easy to see that, if $\mu$ is a lower Ahlfors
$\beta$-regular measure in $\Omega$ then there exists
$\epsilon_0>0$ such that $\mu(\Omega_{r,\epsilon}^t)>0$ for every
$\epsilon\leq
\epsilon_0$.\\

 Now we proceed to state our main results.

\medskip

\begin{theorem}\label{rara}
Let $0\leq\alpha<\beta$ and let $p$ be a log-H\"{o}lder continuous
exponent such that $1<p_*\le p(x)\le p^*<\beta/\alpha$. Let $q$
and $s$ be respectively defined by $1/q(x)=1/p(x)-\alpha/\beta$
and $s(x)=1+q(x)/(p(x))'$. Let $\mu$ be a lower Ahlfors
$\beta$-regular measure in $\Omega$ and $r\in(1,s^*)$. If $w$ is
a weight such that $w(\cdot)^{q(\cdot)} \in A_{r-\delta}(\Omega)$,
for some $\delta\in (0,r-1]$ and
$\int_{\Omega-\Omega_{r,\epsilon_{0}}^s}w(x)^{-(p(x))'}\,
d\mu(x)<\infty$ for some $\epsilon_{0} > 0$, then there exists a
positive constant $C=C(\epsilon)$ such that
\begin{equation*}
\|w\,M_{\alpha}f\|_{q(.),\Omega_{r,\epsilon}^s}\leq C
\|w\,f\|_{p(.),\Omega}
\end{equation*}
for every function $f$ such that $wf\in L^{p(\cdot)}(\Omega)$ and
for every $\epsilon\leq \epsilon_0$.

\end{theorem}

\bigskip

Notice that if $supp\ f\subset \Omega_{r,\epsilon}^s$ then the
hypothesis $\int_{\Omega-\Omega_{r,\epsilon}^s}w(x)^{-(p(x))'}\,
d\mu(x)<\infty$ in the theorem above can be removed and we obtain

\begin{corollary}\label{acotacion}
Let $\alpha$, $p$, $q$, $s$ and $\mu$ be as in theorem
$\ref{rara}$ and let $w$ be a weight such that
$w(\cdot)^{q(\cdot)} \in A_{r-\delta}(\Omega)$, for some $\delta$
such that $0<\delta\leq r-1$. Then there exists a positive
constant $C$ such that the inequality
\begin{equation*}
\|w\,M_{\alpha}f\|_{q(.),\Omega_{r,\epsilon}^s}\leq C
\|w\,f\|_{p(.),\Omega_{r,\epsilon}^s}
\end{equation*}
holds for every function $f$ such that $wf\in
L^{p(\cdot)}(\Omega_{r,\epsilon}^s)$ and $supp\ f\subset
\Omega_{r,\epsilon}^s$.
\end{corollary}

\bigskip

From the fact that $\Omega_{s_*-\delta,\epsilon}^s=\Omega$
we immediately obtain the following result .\\

\begin{corollary}\label{acotacion2}
Let $\alpha$, $p$, $q$, $s$ and $\mu$ be as in theorem
$\ref{rara}$. If $w$ is a weight such that $w(\cdot)^{q(\cdot)}
\in A_{s_*-\delta}(\Omega)$, for some $\delta$ such that
$0<\delta<s_*-1$ then there exists a positive constant $C$ such
that the inequality
\begin{equation*}
\|w\,M_{\alpha}f\|_{q(.),\Omega}\leq C \|w\,f\|_{p(.),\Omega}
\end{equation*}
holds for every function $f$ such that $wf\in
L^{p(\cdot)}(\Omega)$.
\end{corollary}

\bigskip

\bigskip

\begin{corollary}\label{coro1}
Let $\beta=n$, and $\alpha$, $p$, $q$ and $s$ be as in theorem
$\ref{rara}$, and $\mu$ be the Lebesgue measure in $\zR^n$. Let
$x_0\in {\Omega}$ and $w(x)=|x-x_0|^{\eta}$, $\eta\in \zR$. If
$-n/q(x_0)< \eta < n/(p(x_0))'$ then there exists a positive
constant $C$ such that the inequality
\begin{equation*}
\|w\,M_{\alpha}f\|_{q(.),\Omega}\leq C \|w\,f\|_{p(.),\Omega}
\end{equation*}
holds for every function $f$ such that $wf\in
L^{p(\cdot)}(\Omega)$.
\end{corollary}

\bigskip

\begin{remark}
The case $\alpha=0$ in the corollary above was proved in
\cite{KS}. The authors also show that the range of $\eta$ is sharp
 by proving that the reciprocal result is true whenever $x_0\in
\Omega$. Following similar arguments, an analogous result can be
obtained for the case $\alpha>0$.
\end{remark}

\bigskip

An application of corollaries $\ref{acotacion2}$ and $\ref{coro1}$
allows us to obtain two results concerning one-weighted-type
inequalities for the fractional integral operator defined in
$(\ref{intfrac})$.

\bigskip

 Let $0<\alpha<\beta$. Let $p$ be an exponent such that
 $1<p_*<p(x)<p^*<\infty$ for every $x\in\Omega$ and and let $q$ be the function defined
 by $1/q(x)=1/p(x)-\alpha/\beta$.

If $0<\epsilon<\min\{\alpha,
\beta-\alpha, \beta/q^*,\beta(1/p^*-1/q_*)\}$, let
$q_{\epsilon}^{+}$, $q_{\epsilon}^{-}$, $s_{\epsilon}^+$ and
$s_{\epsilon}^-$ be the functions defined by
\begin{equation*}
\frac{1}{q_{\epsilon}^{+}(x)}=\frac{1}{p(x)}-\frac{\alpha+\epsilon}{\beta},
\hspace{1cm}
\frac{1}{q_{\epsilon}^{-}(x)}=\frac{1}{p(x)}-\frac{\alpha-\epsilon}{\beta},\\
\end{equation*}
and
\begin{equation*}
s_{\epsilon}^+(x)=1+\frac{q_{\epsilon}^{+}(x)}{(p(x))'}
 \hspace{2cm}
s_{\epsilon}^-(x)=1+\frac{q_{\epsilon}^{-}(x)}{(p(x))'}.\\
\end{equation*}

\begin{theorem}\label{ialfa}
Let $0<\alpha<\beta$ and $\mu$ be be a lower Ahlfors
$\beta$-regular measure in $\Omega$. Let $p$ be a log-H\"{o}lder
continuous exponent such that $1<p_*\le p(x)\le p^*<\beta/\alpha$
and $q$ be defined by $1/q(x)=1/p(x)-\alpha/\beta$ . Let $w$ be a
weight such that $w^{q_{\epsilon}^{+}}\in
A_{(s_{\epsilon}^+)_{*}}$ and $w^{q_{\epsilon}^{-}}\in
A_{(s_{\epsilon}^-)_{*}}$. Then there exists a positive constant
$C$ such that
\begin{equation*}
\|wI_{\alpha}f\|_{q(\cdot), \Omega}\leq C \|wf\|_{p(\cdot),
\Omega}.
\end{equation*}
\end{theorem}
 \bigskip
   In the classical Lebesgue context and for a particular class
  of weights in $A_1$, we prove that certain variable powers of such weights
  also remain in that class. This property allows us to obtain weights for which
  the boundedness of the fractional integral holds. In order to make the results more
  precise we introduce those special weights.
  \bigskip

  \bigskip

  Let $Q_0$ be a cube in $R^n$ and let $\mu$ be the standard Lebesgue measure.
 We shall be interested in those weights $w$ belonging to the $A_1(Q_0)$ for which the
 following properties hold:\\
 \noindent i) For almost every $x\in Q_0$, $w(x)\geq1$.\\
 \noindent ii) The weight $w$ has $m$ singularities $x_1,\, x_2,\,...,\,x_m$ in $Q_0$. \\
 \noindent iii) There exist two positive numbers $\theta$ and $r$
 such that $w(x) \leq |x-x_i|^{-\theta}$, for almost every $x\in
 Q(x_i,r)\bigcap \bar{Q_0} $ and for each $i = 1,2,...,m$.
\bigskip

 \begin{theorem}\label{reverse}
 Let ($Q_0$,$\mu$) be the measure space consisting of the cube $Q_0$ in $R^n$ and of the classical
  Lebesgue measure $\mu$. Let $w$ be a weight in the $A_1(Q_0)$ class that satisfies the properties stated above.
 \\Then there exists a positive number $\delta$ such that $w^{\alpha}$ also belongs to the $A_1(Q_0)$ class for
 every function $\alpha$ satisfying both a log-H\"{o}lder condition and the inequality $1\leq\alpha(x)\leq1+\delta$
 for almost every point x in $Q_0$.
 \end{theorem}
 \bigskip

\begin{corollary}\label{samko}
Let $\alpha$, $p$ and $q$ be as in corollary $\ref{coro1}$, and
$\mu$ be the Lebesgue measure in a cube $Q_{0}$. For a pair of
weights $w_1$ and $w_2$ in the $A_1(Q_{0})$ class that satisfy the
hypotheses of theorem $\ref{reverse}$ and for $\epsilon>0$ small
enough, let $w$ be the weight defined by
$w=w_1^{1/{q_{\epsilon}^{-}}}w_2^{(1/{q_{\epsilon}^{-}})(1-(s_{\epsilon}^{-})_*)}$,
then there exists a positive constant $C$ such that the inequality
\begin{equation*}
\|w\,I_{\alpha}f\|_{q(.),Q_{0}}\leq C \|w\,f\|_{p(.),Q_{0}}
\end{equation*}
holds for every function $f$ such that $wf\in
L^{p(\cdot)}(Q_{0})$.
\end{corollary}
\begin{remark}When $w$ is a product of a finite number of  power weights, the corollary above can be proved using similar techniques,
for $Q_{0}$ replaced by a more general bounded set $\Omega$.
Particularly, when $w$ is a power weight, this result was proved
by Samko (\cite{S}) but for a variable index $\alpha$.
\end{remark}

\bigskip

The structure of this paper goes on as follows. Section 2 contains
certain types of inequalities frequently used in the variable
context. A pointwise estimate relating both operators $M$ and
$M_{\alpha}$ is also shown. Section 3 is devoted to proving our
main results. Finally, an example of weights in a bounded subset equipped with a non
doubling measure is also given.

\section{Preliminary results}

The following result is a version for cubes of Lemma 3.6 in
\cite{HHP} relating lower Ahlfors regularity and log-H\"{o}lder
continuity of the exponent and it is essential to prove theorem
$\ref{cinco}$. We omit its proof since it is similar to the one
for balls.
\begin{lemma}
Let $\mu$ be a lower Ahlfors $\beta$-regular measure and let $p$
be a log-H\"{o}lder continuous exponent. Then, there exists a
positive constant $C$ such that
\begin{equation}\label{tres}
\mu(Q)^{(p_{Q})_*-(p_Q)^*}\leq C
\end{equation}
for every cube $Q$ centered at $\Omega$.
\end{lemma}

\bigskip

\begin{theorem}\label{cinco}
Let $\mu$ be a lower Ahlfors $\beta$-regular measure in a bounded
$\mu$-measurable set $\Omega$. Let $0<t<t^*$ be a log-H\"{o}lder
continuous function in $\Omega$ and $\epsilon>0$. Then there
exists a positive constant $C=C(\epsilon)$ such that the
inequality
\begin{equation}
Mf(x)^{t(x)}\leq C\left(1+M(|f|(\cdot)^{t(\cdot)})(x)\right)
\end{equation}
holds for almost every $x\in \Omega_{1,\epsilon}^{t}$ and for
every function $f$ such that
$\|f\|_{t(\cdot),\Omega_{1,\epsilon}^{t}}\leq 1$ and
$\int_{\Omega-\Omega_{1,\epsilon}^{t}}|f|\, d\mu\leq 1$, where
$\Omega_{1,\epsilon}^{t}$ is defined as in $(\ref{omega})$.
\end{theorem}

\bigskip

\begin{remark}
The fact that $t$ is allowed to take positive values is essential
in the proof of theorem $\ref{rara}$. If $t_*>1$ and $\mu$ is the
standard Lebesgue measure the inequality above is proved in
\cite{KS}. Inequalities of this type were originally proved by
\cite{D}.
\end{remark}

\bigskip

\Proof  \ : Let $x$ be a fixed point in $\Omega_{1,\epsilon}^t$.
If $Q$ is any cube centered at $x$, it is enough to prove that the
inequality

\begin{equation}\label{cuatro}
\left(\frac{1}{\mu(Q)}\int_{Q}|f(y)|\, d\mu(y)\right)^{{t}(x)}\leq
C \left(1+\frac{1}{\mu(Q)}\int_{Q}|f(y)|^{{t}(y)}\, d\mu(y)\right)
\end{equation}
holds for some positive constant $C$.

Let us first assume that $\mu(Q)\ge 1/2$. Note that
\begin{eqnarray*}
\int_Q |f(y)|\, d\mu(y)=\int_{Q\cap \Omega_{1,\epsilon}^t}
|f(y)|\, d\mu(y)+\int_{Q\cap(\Omega-\Omega_{1,\epsilon}^t)}
|f(y)|\, d\mu(y).
\end{eqnarray*}
The second term in the inequality above is bounded by $1$ by
hypothesis. For the first we take into account that in
$Q\cap\Omega_{1,\epsilon}^t$ we have that $t(x)>1$. Then, from the
fact that

\begin{equation*}
\int_{Q\cap\Omega_{1,\epsilon}^t}
\left(\frac{1}{\mu(Q)+1}\right)^{(t(x))'}\, d\mu(x)\leq \int_Q
\frac{1}{\mu(Q)+1}\, d\mu(x)< 1,
\end{equation*}
we have  that $\|\chi_Q\|_{t'(.),{Q\cap\Omega_{1,\epsilon}^t}}\leq
\mu(Q)+1$. Then, by the generalized H\"{o}lder's inequality
($\ref{holder}$), the remark above and the hypotheses we obtain
\begin{eqnarray*}
\left(\frac{1}{\mu(Q)}\int_{Q\cap\Omega_{1,\epsilon}^t} |f(x)|\,
d\mu(x)\right)^{{t}(x)} &\leq&
\frac{1}{\mu(Q)^{t(x)}}\|f\|_{t(.),{Q\cap\Omega_{1,\epsilon}^t}}^{{t}(x)}\|\chi_Q\|_{t'(.),
{Q\cap\Omega_{1,\epsilon}^t}}^{{t}(x)}\\
&\leq&\left(1+\frac{1}{\mu(Q)}\right)^{{t}(x)}\\
&\leq&C.
\end{eqnarray*}

Now we assume that $\mu(Q)<1/2$ and $l(Q)>C\epsilon$ where $C$ is
a constant depending on the dimension. Then, from the definition
of $\mu$ a constant $C$ depending on $\epsilon$ and the dimension
can be found so that $\mu(Q)\ge C$. Then we proceed as in the
 case above to obtain the result.

If $\mu(Q)<1/2$ and $l(Q)<C\epsilon$ then, it is easy to check
that $(\Omega-\Omega_1) \bigcap Q=\emptyset$. Then $t(x)>1$ in
$Q$.

If $ t_{Q}=\min_{y\in Q}  t(y)$, by applying H\"{o}lder's
inequality we obtain

\begin{eqnarray}\label{uno}
\left(\frac{1}{\mu(Q)}\int_{Q}|f(y)|\, d\mu(y)\right)^{t(x)}&\leq&
\frac{1}{\mu(Q)^{t(x)/t_{Q}}}\left(\int_{Q}|f(y)|^{t_{Q}}\,
d\mu(y)\right)^{\frac{ t(x)}{ t_{Q}}}.\\\nonumber
\end{eqnarray}
Since
\begin{eqnarray*}
\int_{Q}|f(y)|^{ t_{Q}}\, d\mu(y)&= &\int_{Q\cap\{|f|\leq
1\}}|f(y)|^{ t_{Q}}\, d\mu(y)+\int_{Q\cap\{|f|> 1\}}|f(y)|^{
t_{Q}}\,
d\mu(y)\\
&\leq& 2\left(\mu(Q)+\frac{1}{2}\int_{Q}|f(y)|^{ t(y)}\,
d\mu(y)\right),\\
\end{eqnarray*}
and as the expression in brackets is less than 1, from
$(\ref{uno})$ we get
\begin{eqnarray*}
\left(\frac{1}{\mu(Q)}\int_{Q}|f(y)|\, d\mu(y)\right)^{
t(x)}&\leq& C \mu(Q)^{1- t(x)/
t_{Q}}\left(1+\frac{1}{\mu(Q)}\int_{Q}|f(y)|^{ t(y)}\,
d\mu(y)\right).\\
\end{eqnarray*}
But this is ($\ref{cuatro}$)  because of (\ref{tres}). $\square$

\vspace{1cm}

The following lemma gives a pointwise estimation relating both
operators $M$ and $M_{\alpha}$ and it proves to be essential to
obtain our main results.

\bigskip

\begin{lemma}\label{malfapuntual}
Let $\mu$ be a lower Ahlfors $\beta$-regular measure in $\Omega$.
Let $0<\alpha<\beta$ and $p$ be an exponent such that $1<p_*\le
p(x)\le p^*<\beta/\alpha$. Let $q$ be defined by
$1/q(x)=1/p(x)-\alpha/\beta$. If $s(x)=1+q(x)/p'(x)$, then the
following inequality
\begin{equation*}
M_{\alpha}\left(f/w\right)(x)\leq
\left(M(|f|^{p(\cdot)/s(\cdot)}w^{-q(\cdot)/s(\cdot)})(x)\right)^{s(x)/q(x)}
\left(\int _{Q}|f|(y)^{p(y)}\, d\mu(y)\right)^{\alpha/\beta}
\end{equation*}
holds for every function $f$ and for every weight $w$.
\end{lemma}

\bigskip

\Proof : Let $f$ be a non negative function and let $g$ be the
function defined by $g^{s}=f^{p}w^{-q}$. Since
$f/w=g^{s/p}w^{q/p-1}=g^{1-\alpha/\beta}g^{s/p+\alpha/\beta-1}w^{\alpha
q/\beta}$ then, by H\"{o}lder's inequality, we have that
\begin{eqnarray*}
\frac{1}{\mu(Q)^{1-\alpha /\beta}}\int_{Q}\frac{f}{w}\, d\mu&\leq&
 \frac{1}{\mu(Q)^{1-\alpha/\beta}}\int_{Q}g^{s/p}w^{q/p-1}\,
d\mu\\
&\leq& \left(\frac{1}{\mu(Q)}\int_{Q}g\,
d\mu\right)^{1-\alpha/\beta}\left(\int_{Q}g^{(s/p+\alpha
/\beta-1)(\beta/\alpha)}w^{q}\,
d\mu\right)^{\alpha/\beta}\\
\end{eqnarray*}
Since $s/q=1-\alpha/\beta$ and
$(s/p+\alpha/\beta-1)\beta/\alpha=s$ the last expression is
bounded by
\begin{eqnarray*}
&&(Mg(x))^{s(x)/q(x)}\left(\int_{Q}g^{s}w^{q}\,
d\mu\right)^{\alpha/\beta}\\
&& \hspace{3cm} \leq (Mg(x))^{s(x)/q(x)}\left(\int_{Q}f^{p}\,
d\mu\right)^{\alpha/\beta}.\ \square\\
\end{eqnarray*}

\medskip
The following result gives a Welland type inequality in the
context of lower Ahlfors measures. The proof in the euclidean
setting is given in \cite{W}. For more general measures see, for
instance, \cite{GCM}
\begin{lemma}\label{welland}
Let $0<\alpha<\beta$ and $0<\epsilon<\min\{\alpha, \
\beta-\alpha\}$. Then the inequality
\begin{equation}\label{welland1}
|I_{\alpha}f(x)|\leq
C\left(M_{\alpha+\epsilon}f(x)M_{\alpha-\epsilon}f(x)\right)^{1/2}
\end{equation}
holds.
\end{lemma}

\Proof \: Let
$s$ be a positive number. We split $I_{\alpha}$ as follows
\begin{eqnarray}\label{welland2}
I_{\alpha}f(x)&=&\int_{|x-y|<s}\frac{|f(y)||x-y|^{\alpha}}{\mu(Q(x,2|x-y|))}\,
d\mu(y)+\int_{|x-y|\ge
s}\frac{|f(y)||x-y|^{\alpha}}{\mu(Q(x,2|x-y|))}\, d\mu(y)\nonumber\\
&=& I+II.
\end{eqnarray}
By using the property of the measure $\mu$, for the first term we
have
\begin{eqnarray*}
I&=&\sum_{k=0}^{\infty}\int_{2^{-k-1}s\leq|x-y|<2^{-k}s}\frac{|f(y)||x-y|^{\alpha}}{\mu(Q(x,2|x-y|))}\,
d\mu(y)\\
&\leq&
s^{\alpha}\sum_{k=0}^{\infty}\frac{2^{-k\alpha}}{\mu(Q(x,2^{-k}s))}\int_{|x-y|<2^{-k}s}|f(y)|\,
d\mu(y)\\
&\leq&
s^{\alpha}\sum_{k=0}^{\infty}\frac{2^{-k\alpha}}{\mu(Q(x,2^{-k}s))^{(\alpha-\epsilon)/\beta}}\frac{1}{
\mu(Q(x,2^{-k}s))^{1-(\alpha-\epsilon)/\beta}}\int_{|x-y|<2^{-k}s}|f(y)|\,
d\mu(y)\\
&\leq&C
s^{\epsilon} M_{\alpha-\epsilon}f(x)\sum_{k=0}^{\infty}{2^{-k\epsilon}}\\
&\leq&C
s^{\epsilon} M_{\alpha-\epsilon}f(x).\\
\end{eqnarray*}

For the second term we have
\begin{eqnarray*}
II&=&\sum_{k=0}^{\infty}\int_{2^{k}s\leq|x-y|<2^{k+1}s}\frac{|f(y)||x-y|^{\alpha}}{\mu(Q(x,2|x-y|))}\,
d\mu(y)\\
&\leq&
s^{\alpha}\sum_{k=0}^{\infty}\frac{2^{(k+1)\alpha}}{\mu(Q(x,2^{k+1}s))}\int_{|x-y|<2^{k+1}s}|f(y)|\,
d\mu(y)\\
&\leq&
s^{\alpha}\sum_{k=0}^{\infty}\frac{2^{(k+1)\alpha}}{\mu(Q(x,2^{k+1}s))^{(\alpha+\epsilon)/\beta}}\frac{1}{
\mu(Q(x,2^{k+1}s))^{1-(\alpha+\epsilon)/\beta}}\int_{|x-y|<2^{k+1}s}|f(y)|\,
d\mu(y)\\
&\leq&C
s^{-\epsilon} M_{\alpha+\epsilon}f(x)\sum_{k=0}^{\infty}{2^{-k\epsilon}}\\
&\leq&C
s^{-\epsilon} M_{\alpha+\epsilon}f(x).\\
\end{eqnarray*}

Combining both estimates from $(\ref{welland2})$ we obtain
\begin{eqnarray*}
I_{\alpha}f(x)&\leq&C\left(s^{\epsilon}
M_{\alpha-\epsilon}f(x)+s^{-\epsilon}
M_{\alpha+\epsilon}f(x)\right)
\end{eqnarray*}
and thus, by minimizing the expression in brackets in the
inequality above as a function of $s$ we obtain the desired
result. $\square$

\section{Proofs of the main results}

\vspace{0.5cm}

We now proceed to prove our main results.\\

\Proof \  {\it of theorem $\ref{rara}$}: It is enough to prove
that the inequality
\begin{equation*}
\|w\,M_{\alpha}(f/w)\|_{q(.),\Omega_{r,\epsilon}^s}\leq C
\|f\|_{p(.),\Omega}
\end{equation*}
holds for every function $f$ such that $\|f\|_{p(.),\Omega}\leq
C$. But, by lemma $\ref{malfapuntual}$, it is enough to prove that
\begin{equation*}
\|w\,\left(M(|f|^{p/s}w^{-q/s})\right)^{s/q}\|_{q(.),\Omega_{r,\epsilon}^s}\leq
C \|f\|_{p(.),\Omega}
\end{equation*}
which is equivalent to see that the inequality
$\rho_q(\chi_{\Omega_{r,\epsilon}^s}w\,\left(M(|f|^{p/s}w^{-q/s})\right)^{s/q})\leq
C$ holds whenever $\|f\|_{p(.),\Omega}\leq C$.

\bigskip

Let $\tilde{s}(x)=s(x)/r$. Since $w^q \in A_{r}(\Omega)$,
H\"{o}lder's inequality implies that
\begin{eqnarray*}
\rho_{\tilde{s}}(\chi_{\Omega_{r,\epsilon}^s}|f|^{p/s}w^{-q/s})&=&\int_{\Omega_{r,\epsilon}^s}|f|^{p/r}w^{-q/r}\,
d\mu\\
&\leq& C\left(\int_{\Omega}|f|(x)^{p(x)}\,
d\mu(x)\right)^{1/r}\left(\int_{\Omega}w(x)^{-q(x)/(r -1)}\, d\mu(x)\right)^{1/r'}\\
&\leq&C.
\end{eqnarray*}
Thus
$\||f|^{p/s}w^{-q/s}\|_{\tilde{s}(.),\Omega_{r,\epsilon}^s}\leq
C$. On the other hand, the hypothesis on the weight gives
\begin{eqnarray*}
\int_{\Omega-\Omega_{r,\epsilon}^s}|f|^{p/s}w^{-q/s}\, d\mu(x)
&\leq&C
\||f|^{p/s}\|_{s(\cdot),\Omega-\Omega_{r,\epsilon}^s}\|w^{-q/s}\|_{(s(\cdot))',\Omega-\Omega_{r,\epsilon}^s}\\
&\leq&C.
\end{eqnarray*}
Thus theorem $\ref{cinco}$ can be applied by choosing $t(x)=\tilde
s(x)$ and taking into account that
$\Omega_{r,\epsilon}^s=\Omega_{1,\epsilon}^{\tilde s}$. Then we
get
\begin{eqnarray*}
\rho_q(\chi_{\Omega_{r,\epsilon}^s}w\,\left(M(|f|^{p/s}w^{-q/s})\right)^{s/q})
&=&
\int_{\Omega_{r,\epsilon}^s}\left(M\left(|f|^{p/s}w^{-q/s}\right)(x)\right)^{s(x)}w(x)^{q(x)}\,
d\mu(x)\\
&=&\int_{\Omega_{r,\epsilon}^s}\left(\left(M\left(|f|^{p/s}w^{-q/s}\right)(x)\right)^{\tilde{s}(x)}\right)^{r}w(x)^{q(x)}\,
d\mu(x)\\
&\leq&C\int_{\Omega_{r,\epsilon}^s}\left(1+M\left(\left(|f|^{p/s}w^{-q/s}\right)^{\tilde{s}(\cdot)}\right)(x)\right)^{r}w(x)^{q(x)}\,
d\mu(x)\\
&\leq&C+C\int_{\Omega}\left(M\left(\left(|f|^{p/s}w^{-q/s}\right)^{\tilde{s}(\cdot)}\right)(x)\right)^{r}w(x)^{q(x)}\,
d\mu(x).\\
\end{eqnarray*}

From the fact that $w^q\in A_{r-\delta}$, Marcinkiewicz
interpolation theorem can be applied in order to obtain the
boundedness of the maximal operator $M$ in $L^r(w^q)$. Then, from
the estimation above we get that
\begin{eqnarray*}
\rho_q(\chi_{\Omega_{r,\epsilon}^s}w\,\left(M(|f|^{p/s}w^{-q/s})\right)^{s/q})&\leq&C
+C\int_{\Omega}|f|^{p(x)}\,
d\mu(x)\\
&\leq& C.\ \square
\end{eqnarray*}

\vspace{1cm}

\Proof \ {\it of corollary $\ref{coro1}$}:  Without loss of generality we may assume
that $\|f\|_{p(.),\Omega}=1$. Thus, we have to prove
that $\|w\,M_{\alpha}(f/w)\|_{q(.),\Omega}\leq C$.

Since $p$ is a log-H\"{o}lder continuous exponent it is easy to
check that so is $q$ and $w(x)^{q(x)}\sim w(x)^{q(x_0)}$. In fact,
this can be obtained from the property of continuity of $q$ if
$|x-x_0|\le 1/2$. Otherwise the statement is immediate because the
boundedness of $\Omega$.

Since $n/q(x_0)<\beta<n/(p(x_0))'$ then  $w^{q}\in
A_{s(x_0)}(\zR^n)$ and there exists a positive number $\eta$ such
that $w^{q}\in A_{s(x_{0})-\eta}(\zR^n)$; particularly $w^{q}\in
A_{s(x_{0})-\eta}(\Omega)$. By virtue of the continuity of $p$ and
$s$ two positive numbers $\delta$ and $\epsilon$ can be chosen
satisfying $Q(x_0, \delta)\subset
\Omega_{s(x_{0})-\eta/2,\epsilon}^s$ and $\beta(p(x))'<n$ for
$x\in Q(x_0, \delta)$. Thus, we have
\begin{eqnarray*}
\|w\,M_{\alpha}(f/w)\|_{q(.),\Omega}&\leq&
\|w\,M_{\alpha}(f/w)\|_{q(.),Q(x_{0},
\delta)}+\|\chi_{\Omega\backslash Q(x_{0},
\delta)}w\,M_{\alpha}(f/w)\|_{q(.),\Omega}\\
&\leq&\|w\,M_{\alpha}(f/w)\|_{q(.),\Omega_{s(x_{0})-\eta/2,\epsilon}^s}+
\|\chi_{\Omega\backslash Q(x_{0},
\delta)}w\,M_{\alpha}(f/w)\|_{q(.),\Omega}.\\
\end{eqnarray*}
It is easy to see that the hypotheses on the weight are satisfied
as $\Omega$ is bounded and $\Omega\backslash\Omega_{
s(x_{0})-\eta/2,\epsilon}^s\subset \{x\in \Omega:
|x-x_0|>C\delta\}$. Thus theorem $\ref{rara}$ can be applied to
estimate the first term.

If $\beta\leq 0$, $w$ is bounded below for $x\in \Omega\backslash
Q(x_{0}, \delta)$ and its reciprocal is bounded above in $\Omega$
and the boundedness of the second term follows by virtue of the
unweighted norm inequality for $M_{\alpha}$. On the other hand, if
$\beta>0$, let $Q$ be a cube such that $l(Q)\leq \delta/2$. Since
$|y-x_0|\ge \delta/2$ when $y\in Q$ then
\begin{eqnarray*}
\frac{|x-x_0|^{\beta}}{|Q|^{1-\alpha/n}}\int_Q
\frac{f(y)}{|y-x_0|^{\beta}}\, dy\leq
C(diam\Omega)^{\beta}M_{\alpha}f(x).
\end{eqnarray*}
On the other hand, if $l(Q)\ge \delta/2$ we have
\begin{eqnarray*}
\frac{|x-x_0|^{\beta}}{|Q|^{1-\alpha/n}}\int_Q
\frac{f(y)}{|y-x_0|^{\beta}}\, dy&\leq& (diam\Omega)^{\beta}(I+II)
\end{eqnarray*}
where
\begin{equation*}
I=\frac{1}{|Q|^{1-\alpha/n}}\int_{Q\cap \{|y-x_0|\leq \delta/2\}}
\frac{f(y)}{|y-x_0|^{\beta}}\, dy
\end{equation*}
and
\begin{equation*}
II=\frac{1}{|Q|^{1-\alpha/n}}\int_{Q\cap \{|y-x_0|\geq \delta/2\}}
\frac{f(y)}{|y-x_0|^{\beta}}\, dy.
\end{equation*}
It is easy to see that $II\leq C_{\delta}M_{\alpha}f(x)$. Thus we
proceed to estimate $I$. Let $\bar{p}=(p_{Q(x_0,\delta/2)})_*$, by
applying H\"older inequality and taking into account that $\beta
\bar{p}' <n$, we obtain that
\begin{eqnarray*}
I&\leq&\frac{1}{|Q|^{1-\alpha/n}}\left(\int_{Q\cap\{|y-x_0|\leq
\delta/2\}}
|f|^{\bar{p}}\right)^{1/\bar{p}}\left(\int_{\{|y-x_0|\leq
\delta/2\}} |y-x_0|^{-\beta\bar{p}'}\right)^{1/\bar{p}'}\\
&\leq&C \left(\frac{1}{|Q|}\int_{Q\cap \{|y-x_0|\leq \delta/2\}}
|f|^{\bar{p}}\right)^{1/\bar{p}}\\
&\leq&C \left(\frac{1}{|Q|}\int_{Q\cap\{|f|\leq 1\}}
|f|^{\bar{p}}+C_{\delta}\int_{\{|f|\ge 1\}\cap \{\{|y-x_0|\leq
\delta/2\}\}}
|f|^{\bar{p}}\right)^{1/\bar{p}}\\
&\leq& \left(C+C_{\delta}\int_{\{|f|\ge 1\}\cap \{\{|y-x_0|\leq
\delta/2\}\}}
|f|^{{p(y)}}\, dy\right)^{1/\bar{p}}\\
&\leq& C
\end{eqnarray*}
where in the last inequality we have used that
$\|f\|_{p,\Omega}=1$. Thus we have the pointwise inequality
\begin{eqnarray*}
w(x)M_{\alpha}(f/w)(x)\leq C+ CM_{\alpha}f(x),
\end{eqnarray*}
which allows us to obtain the desired result by using the
unweighted classical boundedness of $M_{\alpha}$ and the fact that
$\Omega$ is bounded.

\bigskip

\Proof \ {\it of theorem $\ref{ialfa}$}: It is enough to prove
that the inequality
\begin{equation*}
\|w\,I_{\alpha}(f/w)\|_{q(.),\Omega}\leq C \|f\|_{p(.),\Omega}
\end{equation*}
holds for every function $f$ such that $\|f\|_{p(.),\Omega}\leq
C$.

If we define $q^+(x)=2q_{\epsilon}^+/q(x)$ and
$q^-(x)=2q_{\epsilon}^-/q(x)$ then
$\frac{1}{q^+(x)}+\frac{1}{q^-(x)}=1$. By applying Welland's
inequality and Young's inequality  we obtain
\begin{equation}\label{estimacion}
\int_{\Omega}\left|I_{\alpha}\left(\frac{f}{w}\right)\right|^{q}
w^{q}\, d\mu \hspace{10cm}
\end{equation}

\begin{eqnarray*}
&\leq&
C\left(\int_{\Omega}\frac{1}{q^+}M_{\alpha+\epsilon}\left(\frac{f}{w}\right)^{qq^+/2}
w^{qq^+/2}d\mu+
\int_{\Omega}\frac{1}{q^-}M_{\alpha-\epsilon}\left(\frac{f}{w}\right)^{qq^-/2}
w^{qq^-/2}d\mu \right)\\
&\leq&
C\left(\int_{\Omega}M_{\alpha+\epsilon}\left(\frac{f}{w}\right)^{q_{\epsilon}^+}
w^{q_{\epsilon}^+}d\mu+\int_{\Omega}M_{\alpha-\epsilon}\left(\frac{f}{w}\right)^{q_{\epsilon}^-}
w^{q_{\epsilon}^-}d\mu \right).
\end{eqnarray*}
Now the desired inequality follows immediately because of the
hypothesis on the weights and by virtue of corollary
$\ref{acotacion2}$.

\bigskip
 \Proof \ {\it of theorem $\ref{reverse}$}:

 Let us see that $w^{\alpha}\in A_1(Q_0)$. In fact, a positive constant $\tau$
 can be chosen in such a way that, if $\mu(Q)\leq \tau$ then either $Q$ does contain only one
 singularity or it does not contain any at all.\\
  If $Q$ contains no singularity
 it is easy to see that $w\cong C$ and then, from the fact that
 $w\ge1$ we obtain the result. Now let $x_i$ be the only
 singularity contained in $Q$.  Since $1<w(x)\leq |x-x_i|^{-\theta}$
 and $\alpha$ satisfies a log-H\"{o}lder condition, by taking
 logarithms, we obtain
 \begin{eqnarray*}
 0\leq |\alpha(x)-\alpha(x_i)| \log w \leq C
 \end{eqnarray*}
 which allows us to immediately obtain that $w^{\alpha(x)}\cong
 w^{\alpha(x_i)}$ for almost every $x\in Q$. Thus the result
 follows easily whenever $\mu(Q)$ is small enough.\\
 Let us now consider those cubes $Q$ for which $\mu(Q)> \tau$.
 By a well-known property
 of Muckenhoupt classes there exists a positive number $\delta$ such that $w^{1+\delta}\in A_1(Q_0)$.
 Let's see that this number does work. Let $\alpha$ be a function as in the hypothesis. Since
  $1<\alpha(x)<1+\delta$, for almost every $x$ in $Q$, we have
 \begin{eqnarray*}
 \frac{w^{\alpha}(Q)}{\mu(Q)}&\leq&\frac{w^{1+\delta}(Q)}{\mu(Q)}\\
 &\leq& C\left(\frac{w(Q)}{\mu(Q)}\right)^{1+\delta}\\
 &\leq& C w(\Omega)^{1+\delta}\\
 &\leq& C w(x)^{\alpha(x)}.\square
 \end{eqnarray*}

\Proof \ {\it of corollary $\ref{samko}$}: The thesis follows by
observing that $w$ satisfies the hypothesis in theorem
$\ref{ialfa}$ in the context of the measure space ($Q_0$,$\mu$),
where $Q_0$ is a cube and $\mu$ is the Lebesgue measure. As $w_1$
and $w_2$ $\in A_1$,
 there exist two positive numbers $\delta_1$ and $\delta_2$
such that both weights $w_1^{1+\delta_1}$ and $w_2^{1+\delta_2}$ also belong to that class.
We choose $\delta=\min\{\delta_1,\delta_2\}$.
Now let $\epsilon$ be a positive number as defined in theorem $\ref{ialfa}$ and
 let $\alpha_{\epsilon}$ be the function defined by
 $\alpha_{\epsilon}(x)=q_{\epsilon}^{+}(x)/q_{\epsilon}^{-}(x)$.
 If $\epsilon<\delta\beta/((2+\delta)q^*)$, we then have that
 $1<\alpha_\epsilon<1+\delta$. Moreover,
 since $q$ has the log-H\"{o}lder property, it is easy to see that
 so does $\alpha_{\epsilon}$. Thus $w^{q_\epsilon^{-}}\in
 A_{(s_{\epsilon}^{-})_*}$ since  $w_1 w_2^{1-(s_{\epsilon}^{-})_*} \in
 A_{(s_{\epsilon}^{-})_*}$.
 Moreover
  theorem $\ref{reverse}$ can now be applied to conclude that
 $w^{q_{\epsilon}^{+}}=w_1^{\alpha_{\epsilon}}w_2^{\alpha_{\epsilon}(1-(s_{\epsilon}^{-})_*)}$
  $\in A_{(s_{\epsilon}^{-})_*}$. By virtue of the monotonic
 character of the Muckenhoupt classes it is also true that
 $w^{q_{\epsilon}^{+}}$  $\in A_{(s_{\epsilon}^+)_*}$ and we are
 done. $\square$
 \bigskip

\vspace{1cm}
 We finally obtain a family of weights in the  $A_{p(.)}(\Omega)$ class where
 $\Omega$ has been equipped with a measure $\mu$ that fails to
 have the standard doubling property.\\

Let $X_1=\{(x,x), x\in (0,1)\}$, $X_2=(-1,0)^2$ and
$\Omega=X_1\cup X_2$. If $\Omega$ is any cube containing $\Omega$
and $\mu_i$ is the $i-$dimensional Lebesgue measure for $i=1,2$,
let $\mu$ be the measure supported in $\Omega$ and defined by
$\mu=\mu_i$ in $X_i$  for $i=1,2$. It is easy to prove that $\mu$
is lower Ahlfors 2-regular. If $1<p<\infty$, let $w$ be the weight
defined in $\Omega$ by
\begin{equation*}
w(x,y)= \left\{
 \begin{array}{lll}
 x^{\alpha} & \mbox{if $(x,y)\in X_1$}, \\
 |xy|^{\alpha} & \mbox{if $(x,y)\in X_2$},
 \end{array}
 \right.
\end{equation*}
with $-1<\alpha<p-1$. Then $w$ belongs to the class $A_p(\Omega)$.
In fact, if $Q$ is a cube contained in $\Omega$, it might happen
that $Q$ is a proper subset of $X_1$ or of $X_2$, otherwise $Q$
intersects both of them. In the first two cases the statement
follows immediately  from the $A_p$ conditions for the
1-dimensional and 2-dimensional Lebesgue measure respectively.

Let us prove the remaining case. Given $b\leq a<0$, and $l>0$ such
that $0<b+l<a+l$ let $Q=(a,a+l)\times (b,b+l)$. Let $I=(a,0)$,
$J=(b,0)$ and $K=(0,b+l)$ thus $Q\cap\Omega=(I\times J)\cup
\{(x,x), x\in K\}$ and $\mu(Q)=\mu(Q\cap\Omega)=\mu_2(I\times
J)+\sqrt{2}\mu_1{(K)}$. By the definition of $w$ and Tonelli's
theorem we obtain

\begin{eqnarray*}
&&\left(\frac{1}{\mu(Q)}\int_{Q}w\,d\mu\right)\left(\frac{1}{\mu(Q)}\int_{Q}w^{-1/(p-1)}\,
d\mu\right)^{p-1}\hspace{11cm} \\
\end{eqnarray*}
\begin{eqnarray*}
&& \leq\frac{1}{\mu_1(I)^p}\left(\int_{I}|x|^{\alpha}dx
\right)\left(\int_{I}|x|^ {-\alpha/(p-1)}
dx\right)^{p-1}\frac{1}{\mu_1(J)^p} \left(\int_{J}|x|^{\alpha}
dx\right)\left(\int_{J}|x|^
{-\alpha/(p-1)} dx\right)^{p-1}\\
&&
\hspace{0.5cm}+\frac{1}{\mu_1(K)^p}\left(\int_{K}|x|^{\alpha}dx\right)\left(\int_{K}|x|^
{-\alpha/(p-1)}dx\right)^{p-1}\\
&&\hspace{0.5cm}+C\left(\frac{\mu_1(I)\mu_1(J)}{\mu_1(I)\mu_1(J)+\sqrt{2}\mu_1(K)}\right)^{\alpha+1}
\left(\frac{\sqrt{2}\mu_1(K)}{\mu_1(I)\mu_1(J)+\sqrt{2}\mu_1(K)}\right)^{p-\alpha-1}\\
&&\hspace{0.5cm}+C\left(\frac{\mu_1(I)\mu_1(J)}{\mu_1(I)\mu_1(J)+\sqrt{2}\mu_1(K)}\right)^{p-\alpha-1}
\left(\frac{\sqrt{2}\mu_1(K)}{\mu_1(I)\mu_1(J)+\sqrt{2}\mu_1(K)}\right)^{\alpha+1}\\
&&\leq C
\end{eqnarray*}
where we have used the one dimensional $A_p$ inequality for both
the first and the second terms and the range of $\alpha$ for the
boundedness of the last two terms.

\end{document}